# Efficient Designs in Small Blocks for Comparing Consecutive Pairs of Treatments

S. Huda  and  Rahul Mukerjee
Department of Statistics and OR         Indian Institute of Management Calcutta
Faculty of Science, Kuwait University   Joka, Diamond Harbour Road
P.O. Box-5969, Safat-13060, Kuwait      Kolkata 700 104, India

**Abstract**: Optimal block designs in small blocks are explored when the treatments have a natural ordering and interest lies in comparing consecutive pairs of treatments. We first develop an approximate theory which leads to a convenient multiplicative algorithm for obtaining optimal design measures. This, in turn, yields highly efficient exact designs even when the number of blocks is rather small. Moreover, our approach is seen to allow nesting of such efficient exact designs which is an advantage when the resources for the experiment are available possibly in several stages. Illustrative examples are given. Tables of optimal design measures are also provided.

## 1. Introduction and preliminaries

Optimal block designs have received significant attention in the experimental design literature. Commonly, these designs have been studied for inference on complete sets of orthnormal treatment contrasts or contrasts representing treatment versus control comparisons; see e.g., Dey (2010, Ch. 6) for a review and further references. In the present work, we consider experimental situations where the treatments have a natural ordering and interest lies in comparing consecutive pairs of treatments. This happens, for instance, when the treatments are doses of, say, a fertilizer or a drug, or correspond to successive points of time as in time course microarray experiments (Schiffl and Hilgers, 2012). The resulting optimal block design problem opens up new challenges because the contrasts of interest are now quite different from the ones traditionally focused on.

Specifically, suppose there $v$ naturally ordered treatments $1,\ldots,v$, the treatment contrasts of interest being $\tau_i - \tau_{i-1}$, $2 \le i \le v$, where $\tau_1,\ldots,\tau_v$ are the effects of treatments $1,\ldots,v$. Let $L\tau$ be the matrix representation of these contrasts, where $L$ is a $(v-1) \times v$ matrix of full row rank, $\tau = (\tau_1,\ldots,\tau_v)'$, and the prime denotes transposition. Consider a block design $d$, which lays out the $v$ treatments in $b$ blocks each of size $k$ ($< v$). Let $N$ be the treatment versus block incidence matrix of $d$ and $R = \mathrm{diag}(r_1,\ldots,r_v)$, where $r_1,\ldots,r_v$ are the replication numbers of treatments $1,\ldots,v$ in $d$. Then, under the usual fixed effects additive linear model with homoscedastic and uncorrelated observational errors, the information matrix for $\tau$ in $d$ is given by $C = R - k^{-1}NN'$. We consider only designs $d$ which keep $L\tau$ estimable. For any such $d$, the matrix $LCL'$ is nonsingular and one can show that $L =$ $(LL')(LCL')^{-1}LC$. Hence, writing $L\hat{\tau}$ for the best linear unbiased estimator of $L\tau$, one obtains $\text{cov}(L\hat{\tau}) = \sigma^2 W^{-1}$, where $\sigma^2$ is the constant error variance, and

$$W = (LL')^{-1} LCL'(LL')^{-1}. \qquad (1)$$

Given $v$, $k$ and $b$, we intend to explore $d$ so as to keep the average variance of the elements of $L\hat{\tau}$, or equivalently $\text{tr}(W^{-1})$, small. This corresponds to the A-criterion which is appropriate here because interest is focused on the contrasts in $L\tau$ as they stand, and not in their linear functions; cf. Kerr (2012) and Schiffl and Hilgers (2012). In general, however, the task of minimizing $\text{tr}(W^{-1})$ over all possible designs $d$ is extremely involved due to complexity of the underlying combinatoics. The difficulty is accentuated by the fact that the matrix $T = (LL')^{-1}L$ in $W$ does not have a simple form. Moreover, as seen later, in the process one cannot ignore nonbinary blocks which accommodate one or more treatments more than once.

We show how significant progress can be made on this challenging problem via the approximate theory which is seen to provide a useful benchmark and yield exact designs with assured efficiency often over 0.95 or even 0.99. For most practical purposes, it suffices to have designs with such assured high efficiency. Moreover, we demonstrate that this is achievable even in a relatively small number of blocks which is important from the perspective of experimental economy. In addition, our approach is seen to enable nesting of such efficient exact designs which is an advantage when the resources for the experiment are available possibly in several stages.

**2. Approximate theory**

With $v$ treatments, let $\mathcal{B} = \{S_1,...,S_B\}$ be the class of all possible blocks of size $k$, including nonbinary blocks but excluding blocks which include only one treatment and hence carry no information on the treatment contrasts. The blocks $S_1,...,S_B$ in $\mathcal{B}$ play the role of design points in the approximate theory developed here. Suppose treatments $1,...,v$ appear $h_{1j},...,h_{vj}$ times, respectively, in the block $S_j$, $1 \leq j \leq B$, where the nonnegative integers $h_{1j},...,h_{vj}$ add up to $k$. Let $C_j = H_j - k^{-1}h_j h_j'$ be the counterpart of $C$ for the single-block design given by $S_j$, where $H_j = \text{diag}(h_{1j},...,h_{vj})$ and $h_j = (h_{1j},...,h_{vj})'$. With a block design $d$ involving $b$ blocks of size $k$ each, suppose $S_j$ appears $f_j$ times as a block of $d$, $1 \leq j \leq B$, where $f_1,...,f_B$ are nonnegative integers satisfying $f_1 +...+ f_B = b$. Considering the contribution of the individual blocks to the $C$ matrix of $d$, we then get $C = \Sigma_{j=1}^B f_j C_j$. Hence the matrix $W$ in (1) can be expressed as $W = bM(p)$, where

$$M(p) = \Sigma_{j=1}^B p_j V_j, \qquad (2)$$

with $p_j = f_j/b$, $V_j = TC_jT'$, $1 \leq j \leq B$, $p = (p_1,...,p_B)$, and as before, $T = (LL')^{-1}L$.

Thus, given $v$, $k$ and $b$, the design exerts an influence on $W$, or equivalently $M(p)$, only through the proportions $p_1,...,p_B$. We need to find these so as to minimize $\text{tr}[\{M(p)\}^{-1}]$. While the discreteness of $f_1,...,f_B$ induces the same on $p_1,...,p_B$, considerable simplicity is achieved if, for the time being, we treat $p_1,...,p_B$ as nonnegative continuous variables satisfying $p_1+...+p_B=1$. Any such $p = (p_1,...,p_B)$ is called a design measure which assigns masses $p_1,...,p_B$ on the blocks $S_1,...,S_B$, respectively, in $\mathcal{B}$. The design problem now reduces to finding an A-optimal design measure which minimizes $\text{tr}[\{M(p)\}^{-1}]$ over all possible $p$ such that $M(p)$ is nonsingular.

**Lemma 1**. *A design measure $p$ is A-optimal if and only if $M(p)$ is nonsingular and*

$$\text{tr}[\{M(p)\}^{-1}V_j\{M(p)\}^{-1}] \leq \text{tr}[\{M(p)\}^{-1}], \; 1 \leq j \leq B.$$

*Proof.* Let $e_1,...,e_B$ be the $B \times 1$ unit vectors, and $\phi(p) = \text{tr}[\{M(p)\}^{-1}]$. From consideration of directional derivatives in the spirit of Silvey (1980, pp. 18-20), it suffices to show that for any $p$ with nonsingular $M(p)$,

$$\lim_{\varepsilon \to 0}\{\phi((1-\varepsilon)p + \varepsilon e_j) - \phi(p)\}/\varepsilon = \text{tr}[\{M(p)\}^{-1}] - \text{tr}[\{M(p)\}^{-1}V_j\{M(p)\}^{-1}], \; 1 \leq j \leq B. \quad (3)$$

Each $V_j$, being nonnegative definite, can be expressed as $V_j = Q_jQ_j'$ for some matrix $Q_j$. Hence by (2), $M((1-\varepsilon)p + \varepsilon e_j) = (1-\varepsilon)M(p) + \varepsilon Q_jQ_j'$, so that

$$\{M((1-\varepsilon)p + \varepsilon e_j)\}^{-1}$$
$$= (1-\varepsilon)^{-1}[\{M(p)\}^{-1} - \delta\{M(p)\}^{-1}Q_j(I + \delta Q_j'\{M(p)\}^{-1}Q_j)^{-1}Q_j'\{M(p)\}^{-1}],$$

where $\delta = \varepsilon/(1-\varepsilon)$ and $I$ is the identity matrix with as many columns as $Q_j$. The truth of (3) follows from the above since $V_j = Q_jQ_j'$. □

Although the A-optimal design measure is not analytically tractable, Lemma 1 suggests a multiplicative algorithm which can be used conveniently for its numerical determination. The algorithm starts with the uniform measure $p^{(0)} = (1/B,...,1/B)$, and finds $p^{(u)} = (p_1^{(u)},...,p_B^{(u)})$, $u = 1,2,...$ recursively as

$$p_j^{(u)} = p_j^{(u-1)} \frac{\text{tr}[\{M(p^{(u-1)})\}^{-1}V_j\{M(p^{(u-1)})\}^{-1}]}{\text{tr}[\{M(p^{(u-1)})\}^{-1}]}, \qquad 1 \leq j \leq B, \quad (4)$$

till a design measure $p^{(u)}$, satisfying

$$\text{tr}[\{M(p^{(u)})\}^{-1}V_j\{M(p^{(u)})\}^{-1}] - \text{tr}[\{M(p^{(u)})\}^{-1}] \leq t, \; 1 \leq j \leq B, \quad (5)$$

is reached, where *t* is a preassigned small positive quantity. Write $\tilde{p} = (\tilde{p}_1,...,\tilde{p}_B)$ for the terminal design measure meeting (5). Then arguing alone the lines of Silvey (1980, p. 35), $\text{tr}[\{M(\tilde{p})\}^{-1}]$ exceeds the minimum possible $\text{tr}[\{M(p)\}^{-1}]$ by at most *t*. We take $t = 10^{-10}$. So, $\tilde{p}$ represents the *A*-optimal design measure with accuracy as high as up to nine places of decimals. Even at this level of accuracy, the algorithm in (4) and (5) is seen to be quite fast.

In fact, a slight modification, that takes a chance that the *A*-optimal design measure may be supported only on binary blocks often makes the algorithm even faster. Note that a binary block is one where each treatment appears at most once. Given *v* and *k*, consider the subclass of $\mathcal{B}$ as given by the $B_0 = \binom{v}{k}$ binary blocks. One may then first run a modified version of the algorithm in (4) and (5) to find the *A*-optimal design measure among those with support confined to this subclass. There are two possibilities about the resulting design measure:

(i) It meets (5) for $1 \leq j \leq B$. Then it is as good as what the original version of the algorithm would yield.

(ii) It fails to meet (5) for some *j*, $1 \leq j \leq B$. Then it not *A*-optimal among all design measures and one needs to run the algorithm in (4) and (5) as it stands.

The point here is that if (i) arises, which is often the case, then there is further saving in computational time because $B_0$ is much smaller than *B* for $k \geq 3$.

Tables 1-4 show *A*-optimal design measures for $k < v \leq 10$ and *k* = 2, 3, 4 and 5, respectively. Only blocks with positive masses are displayed in the tables and we write $\tilde{\phi} = \text{tr}[\{M(\tilde{p})\}^{-1}]$. All these design measures, barring the ones for *k* =5 and *v* = 7, 9, are supported only over binary blocks and hence the modified version of the algorithm accelerates their derivation. On the other hand, the two exceptional cases vindicate a point mentioned earlier that in general there is no guarantee that binary blocks alone would suffice for our purpose. The contents of Table 1 for *k* = 2 strengthen the findings in Schiffl (2011) and Schiffl and Hilgers (2012) who also considered approximate theory but stopped short of exploring directional derivatives, and reported results for *v* = 3, 4 and 5.

As the tables reveal, unlike in many other design problems, the *A*-optimal design measure in the present context is far from uniform. Hence it provides useful indication about good exact designs – e.g., it suggests that an efficient exact design should include, possibly with repetition, those blocks where it assigns greater masses. A more precise account follows in the next section.

## 3. Efficient exact designs and examples

Given $v$ and $k$, an exact design $d$ in $b$ blocks corresponds to a design measure $p^{exact} = (f_1/b, ..., f_B/b)$, where $f_j$ is the number of times $S_j$ appears as a block of $d$, and $f_1 + ... + f_B = b$. Hence its $A$-efficiency, relative to the $A$-optimal design measure $\tilde{p}$, can be defined as

$$\text{Eff} = \frac{\text{tr}[\{M(\tilde{p})\}^{-1}]}{\text{tr}[\{M(p^{exact})\}^{-1}]} \tag{6}$$

Thus $\tilde{p}$ acts as a useful benchmark for assessing the performance of exact designs.

Having obtained $\tilde{p} = (\tilde{p}_1, ..., \tilde{p}_B)$ as in the last section, a simple rounding off technique often yields an efficient exact design in $b$ blocks. Suppose there exists a constant $c$ such that the quantities $\tilde{f}_j = [c\tilde{p}_j]$, $1 \leq j \leq B$, add up to $b$, where $[x]$ denotes the integer nearest to any nonnegative $x$. Since $\tilde{f}_1, ..., \tilde{f}_B$ are approximately proportional to $\tilde{p}_1, ..., \tilde{p}_B$, one expects intuitively that an exact design, where $S_j$ appears $\tilde{f}_j$ times as a block, $1 \leq j \leq B$, should be highly efficient. As our examples show, this is really the case even for relatively small $b$, where it is typically harder to match $\tilde{p}$ by an exact design.

For certain combinations of $v$, $k$ and $b$, a constant $c$ as envisaged above may not exist. Examples 4-6 illustrate how a little adjustment of the rounding off technique then helps in finding highly efficient exact designs. These examples also reveal another advantage of our procedure, namely, nesting of designs, which is very helpful when the resources for the experiment are available possibly in several stages. We elaborate on this after presenting Example 5.

In all examples, the stated efficiency figures for exact designs are relative to the $A$-optimal design measure $\tilde{p}$; vide (6). These figures are conservative because $\tilde{p}$ is not attainable in an exact setup. Thus the true efficiencies of the exact designs in the examples are even higher and some of them may well be $A$-optimal, among all exact designs with the same $v$, $k$ and $b$, for comparing consecutive treatments.

**Example 1**. Let $v = 6$, $k = 2$, $b = 14$. The $A$-optimal design measure for $v = 6$, $k = 2$ is as shown in Table 1. If we multiply the masses of this design measure by $c = 14.2$, and then round these off to nearest integers, then we get the vector $(2, 1, 0, 0, 0, 2, 1, 0, 0, 2, 1, 0, 2, 1, 2)'$ whose elements add up to $b$ (=14), and correspond to blocks $\{1,2\}, \{1,3\}, \{1,4\}, ...$, respectively Thus we get the exact design consisting of $b=14$ blocks

$\{1,2\}, \{1,2\}, \{1,3\}, \{2,3\}, \{2,3\}, \{2,4\}, \{3,4\}, \{3,4\}, \{3,5\}, \{4,5\}, \{4,5\}, \{4,6\}, \{5,6\}, \{5,6\}$.

This design has $A$-efficiency 0.9650. □

**Example 2**. This example concerns the smallest ($v$, $k$) such that the support of the $A$-optimal design measure includes nonbinary blocks. Let $v$ =7, $k$ =5, $b$=9. As in Example 1, we multiply the masses of the $A$-optimal design measure shown in Table 4 by $c$ =9, and then round these off to nearest integers to get the nonbinary exact design with blocks

$$\{1,2,3,4,5\}, \{1,2,3,4,5\}, \{1,2,3,6,7\}, \{1,2,5,6,7\}, \{2,3,4,5,6\},$$
$$\{3,4,5,6,7\}, \{3,4,5,6,7\}, \{1,2,2,3,4\}, \{4,5,6,6,7\}$$

and $A$-efficiency 0.9992. □

**Example 3**. We now consider a value of $v$ beyond the range of Tables 1-4 and show that still our approach works quite well. Let $v$=12, $k$=4, $b$=11. The $A$-optimal design measure $\tilde{p}$ turns out to be as follows and has $\text{tr}[\{M(\tilde{p})\}^{-1}]$ =70.8503.

| Block | Mass | Block | Mass | Block | Mass | Block | Mass |
|---|---|---|---|---|---|---|---|
| 1 2 3 4 | 0.1223 | 1 2 8 9 | 0.0055 | 3 10 11 12 | 0.0113 | 5 10 11 12 | 0.0176 |
| 1 2 3 7 | 0.0221 | 1 2 10 11 | 0.0071 | 4 5 6 7 | 0.0797 | 6 7 8 9 | 0.0797 |
| 1 2 3 8 | 0.0176 | 1 2 11 12 | 0.0230 | 4 5 8 9 | 0.0027 | 6 10 11 12 | 0.0221 |
| 1 2 3 9 | 0.0110 | 2 3 4 5 | 0.0729 | 4 5 11 12 | 0.0055 | 7 8 9 10 | 0.0708 |
| 1 2 3 10 | 0.0113 | 2 3 11 12 | 0.0071 | 4 10 11 12 | 0.0110 | 7 8 11 12 | 0.0307 |
| 1 2 5 6 | 0.0307 | 3 4 5 6 | 0.0708 | 5 6 7 8 | 0.0723 | 8 9 10 11 | 0.0729 |
|  |  |  |  |  |  | 9 10 11 12 | 0.1223 |

We multiply the masses of this measure by $c$ =12.3, and then round these off to nearest integers to get the exact design with blocks

$$\{1,2,3,4\}, \{1,2,3,4\}, \{2,3,4,5\}, \{3,4,5,6\}, \{4,5,6,7\}, \{5,6,7,8\},$$
$$\{6,7,8,9\}, \{7,8,9,10\}, \{8,9,10,11\}, \{9,10,11,12\}, \{9,10,11,12\},$$

and $A$-efficiency 0.9562. □

**Example 4**. Let $v$ =6, $k$ =3, $b$=11. Starting from the $A$-optimal design measure shown in Table 2, multiplication by $c$ =10.6928 and 10.6929, followed by rounding off, leads to exact designs

$d(10)$:  $\{1,2,3\}, \{1,2,3\}, \{1,2,4\}, \{2,3,4\}, \{2,3,4\},$
$\{3,4,5\}, \{3,4,5\}, \{3,5,6\}, \{4,5,6\}, \{4,5,6\},$

and

$d(12)$:  $\{1,2,3\}, \{1,2,3\}, \{1,2,3\}, \{1,2,4\}, \{2,3,4\}, \{2,3,4\},$
$\{3,4,5\}, \{3,4,5\}, \{3,5,6\}, \{4,5,6\}, \{4,5,6\}, \{4,5,6\},$

with 10 and 12 blocks, and $A$-efficiencies 0.9547 and 0.9785, respectively. No multiplier yields an exact design with 11 blocks. Observe that that $d(12)$ includes all blocks of $d(10)$ and two additional blocks, namely, $\{1,2,3\}$ and $\{4,5,6\}$. Deletion of any of these two blocks from $d(12)$ produces a design $d(11)$ with $b$=11 blocks and having $A$-efficiency 0.9578. □

**Example 5**. Let $v$ =9, $k$ =4, $b$=10. Starting from the $A$-optimal design measure shown in Table 3, multiplication by $c$ =13.2978 and 13.2979, followed by rounding off, leads to exact designs

$$d(9): \quad \{1,2,3,4\}, \{1,2,3,4\}, \{1,2,8,9\}, \{2,3,4,5\}, \{3,4,5,6\},$$
$$\{4,5,6,7\}, \{5,6,7,8\}, \{6,7,8,9\}, \{6,7,8,9\},$$

and

$$d(11): \quad \{1,2,3,4\}, \{1,2,3,4\}, \{1,2,5,6\}, \{1,2,8,9\}, \{2,3,4,5\}, \{3,4,5,6\},$$
$$\{4,5,6,7\}, \{4,5,8,9\}, \{5,6,7,8\}, \{6,7,8,9\}, \{6,7,8,9\},$$

having 9 and 11 blocks, and *A*-efficiencies 0.9902 and 0.9731, respectively. No multiplier yields an exact design with 10 blocks. Deleting from $d(11)$ any of the two blocks $\{1,2,5,6\}$ and $\{4,5,8,9\}$ that are not in $d(9)$, one gets a design $d(10)$ with $b=10$ blocks and having *A*-efficiency 0.9777. □

In Example 4, the designs $d(10)$, $d(11)$ and $d(12)$ are nested in the sense that all blocks of $d(10)$ appear in $d(11)$ and all blocks of $d(11)$ appear in $d(12)$. This is yet another advantage of approximate theory-aided construction of exact designs when the resources are likely to be available in several stages. Thus in the setup of Example 4, if the initially available resources allow 10 blocks, then one may start with $d(10)$. Later, if additional resources are available in two stages allowing one extra block at each stage, then one may first augment $d(10)$ by the block $\{1,2,3\}$ to get $d(11)$, and subsequently augment $d(11)$ by the block $\{4,5,6\}$ to get $d(12)$, retaining high efficiency at each stage. The same advantage accrues also in Example 5, where the designs $d(9)$, $d(10)$ and $d(11)$, all highly efficient, are nested. We conclude with a more elaborate example of nested exact designs.

**Example 6**. Let $v=10$, $k=5$. Suppose the initially available resources allow 10 blocks but there is a possibility of having additional resources in two stages allowing three extra blocks each stage. Thus one needs to plan for efficient exact designs in 10, 13 and 16 blocks. Starting from the *A*-optimal design measure shown in Table 4, multiplication by $c$ =10, 14.1376, 14.1377 and 18, followed by rounding off, leads to exact designs

$$d(10): \quad \{1,2,3,4,5\}, \{1,2,3,4,5\}, \{1,2,3,9,10\}, \{1,2,8,9,10\}, \{2,3,4,5,6\},$$
$$\{3,4,5,6,7\}, \{4,5,6,7,8\}, \{5,6,7,8,9\}, \{6,7,8,9,10\}, \{6,7,8,9,10\},$$

$$d(12): \quad \{1,2,3,4,5\}, \{1,2,3,4,5\}, \{1,2,3,4,5\}, \{1,2,3,9,10\}, \{1,2,8,9,10\}, \{2,3,4,5,6\},$$
$$\{3,4,5,6,7\}, \{4,5,6,7,8\}, \{5,6,7,8,9\}, \{6,7,8,9,10\}, \{6,7,8,9,10\}, \{6,7,8,9,10\},$$

$$d(14): \quad \{1,2,3,4,5\}, \{1,2,3,4,5\}, \{1,2,3,4,5\}, \{1,2,3,9,10\}, \{1,2,8,9,10\},$$
$$\{2,3,4,5,6\}, \{2,3,4,5,6\}, \{3,4,5,6,7\}, \{4,5,6,7,8\}, \{5,6,7,8,9\},$$
$$\{5,6,7,8,9\}, \{6,7,8,9,10\}, \{6,7,8,9,10\}, \{6,7,8,9,10\},$$

$$d(16): \quad \{1,2,3,4,5\}, \{1,2,3,4,5\}, \{1,2,3,4,5\}, \{1,2,3,4,5\}, \{1,2,3,9,10\},$$
$$\{1,2,8,9,10\}, \{2,3,4,5,6\}, \{2,3,4,5,6\}, \{3,4,5,6,7\}, \{4,5,6,7,8\}, \{5,6,7,8,9\},$$
$$\{5,6,7,8,9\}, \{6,7,8,9,10\}, \{6,7,8,9,10\}, \{6,7,8,9,10\}, \{6,7,8,9,10\},$$

with 10, 12, 14 and 16 blocks, and *A*-efficiencies 0.9951, 0.9955, 0.9920 and 0.9942, respectively. No multiplier yields an exact design with 13 blocks. Now, $d(14)$ includes all blocks of $d(12)$ and two additional blocks $\{2,3,4,5,6\}$ and $\{5,6,7,8,9\}$. Deletion of any of these two blocks from $d(14)$ produces a design $d(13)$ with 13 blocks and having $A$-efficiency 0.9911. Clearly, $d(10)$ is nested in $d(13)$ which, in turn, is nested in $d(16)$. The fact that $d(10)$ has slightly higher $A$-efficiencies than $d(13)$ and $d(16)$ is a consequence of discretization – it only means that $d(10)$ matches the $A$-optimal design measure more closely than $d(13)$ and $d(16)$. The key point is that all the three nested designs have very high $A$-efficiencies. □

**Acknowledgement**. The work of SH was supported by a grant from Kuwait University. The work of RM was supported by the J.C. Bose National Fellowship of the Government of India and a grant from the Indian Institute of Management Calcutta.

Table 1. *A-optimal design measures for block size* $k = 2$

| $v=3$ | | $v=7$ | | $v=8$ | | $v=9$ | | $v=10$ | |
|---|---|---|---|---|---|---|---|---|---|
| $\widetilde{\phi} = 7.4641$ | | $\widetilde{\phi} = 62.0195$ | | $\widetilde{\phi} = 83.7805$ | | $\widetilde{\phi} = 108.7860$ | | $\widetilde{\phi} = 137.0352$ | |
| Block | Mass | Block | Mass | Block | Mass | Block | Mass | Block | Mass |
| 1 2 | 0.4226 | 1 2 | 0.1321 | 1 2 | 0.1130 | 1 2 | 0.0989 | 1 2 | 0.0880 |
| 1 3 | 0.1548 | 1 3 | 0.0337 | 1 3 | 0.0284 | 1 3 | 0.0245 | 1 3 | 0.0216 |
| 2 3 | 0.4226 | 1 4 | 0.0175 | 1 4 | 0.0143 | 1 4 | 0.0122 | 1 4 | 0.0106 |
| $v=4$ | | 1 5 | 0.0118 | 1 5 | 0.0093 | 1 5 | 0.0078 | 1 5 | 0.0066 |
| $\widetilde{\phi} = 16.2195$ | | 1 6 | 0.0094 | 1 6 | 0.0071 | 1 6 | 0.0057 | 1 6 | 0.0048 |
| Block | Mass | 1 7 | 0.0083 | 1 7 | 0.0060 | 1 7 | 0.0046 | 1 7 | 0.0038 |
| 1 2 | 0.2706 | 2 3 | 0.1159 | 1 8 | 0.0055 | 1 8 | 0.0040 | 1 8 | 0.0032 |
| 1 3 | 0.0806 | 2 4 | 0.0280 | 2 3 | 0.0991 | 1 9 | 0.0038 | 1 9 | 0.0029 |
| 1 4 | 0.0537 | 2 5 | 0.0150 | 2 4 | 0.0234 | 2 3 | 0.0866 | 1 10 | 0.0027 |
| 2 3 | 0.2439 | 2 6 | 0.0108 | 2 5 | 0.0121 | 2 4 | 0.0201 | 2 3 | 0.0769 |
| 2 4 | 0.0806 | 2 7 | 0.0094 | 2 6 | 0.0082 | 2 5 | 0.0101 | 2 4 | 0.0177 |
| 3 4 | 0.2706 | 3 4 | 0.1135 | 2 7 | 0.0066 | 2 6 | 0.0067 | 2 5 | 0.0087 |
| $v=5$ | | 3 5 | 0.0271 | 2 8 | 0.0060 | 2 7 | 0.0051 | 2 6 | 0.0056 |
| $\widetilde{\phi} = 28.2360$ | | 3 6 | 0.0150 | 3 4 | 0.0968 | 2 8 | 0.0044 | 2 7 | 0.0042 |
| Block | Mass | 3 7 | 0.0118 | 3 5 | 0.0223 | 2 9 | 0.0040 | 2 8 | 0.0034 |
| 1 2 | 0.2000 | 4 5 | 0.1135 | 3 6 | 0.0116 | 3 4 | 0.0845 | 2 9 | 0.0030 |
| 1 3 | 0.0547 | 4 6 | 0.0280 | 3 7 | 0.0082 | 3 5 | 0.0190 | 2 10 | 0.0029 |
| 1 4 | 0.0317 | 4 7 | 0.0175 | 3 8 | 0.0071 | 3 6 | 0.0096 | 3 4 | 0.0750 |
| 1 5 | 0.0251 | 5 6 | 0.1159 | 4 5 | 0.0963 | 3 7 | 0.0064 | 3 5 | 0.0167 |
| 2 3 | 0.1770 | 5 7 | 0.0337 | 4 6 | 0.0223 | 3 8 | 0.0051 | 3 6 | 0.0082 |
| 2 4 | 0.0481 | 6 7 | 0.1321 | 4 7 | 0.0121 | 3 9 | 0.0046 | 3 7 | 0.0053 |
| 2 5 | 0.0317 | | | 4 8 | 0.0093 | 4 5 | 0.0839 | 3 8 | 0.0041 |
| 3 4 | 0.1770 | | | 5 6 | 0.0968 | 4 6 | 0.0188 | 3 9 | 0.0034 |
| 3 5 | 0.0547 | | | 5 7 | 0.0234 | 4 7 | 0.0096 | 3 10 | 0.0032 |
| 4 5 | 0.2000 | | | 5 8 | 0.0143 | 4 8 | 0.0067 | 4 5 | 0.0744 |
| $v=6$ | | | | 6 7 | 0.0991 | 4 9 | 0.0057 | 4 6 | 0.0163 |
| $\widetilde{\phi} = 43.5040$ | | | | 6 8 | 0.0284 | 5 6 | 0.0839 | 4 7 | 0.0080 |
| Block | Mass | | | 7 8 | 0.1130 | 5 7 | 0.0190 | 4 8 | 0.0053 |
| 1 2 | 0.1589 | | | | | 5 8 | 0.0101 | 4 9 | 0.0042 |
| 1 3 | 0.0416 | | | | | 5 9 | 0.0078 | 4 10 | 0.0038 |
| 1 4 | 0.0225 | | | | | 6 7 | 0.0845 | 5 6 | 0.0744 |
| 1 5 | 0.0161 | | | | | 6 8 | 0.0201 | 5 7 | 0.0163 |
| 1 6 | 0.0138 | | | | | 6 9 | 0.0122 | 5 8 | 0.0082 |
| 2 3 | 0.1399 | | | | | 7 8 | 0.0866 | 5 9 | 0.0056 |
| 2 4 | 0.0352 | | | | | 7 9 | 0.0245 | 5 10 | 0.0048 |
| 2 5 | 0.0202 | | | | | 8 9 | 0.0989 | 6 7 | 0.0744 |
| 2 6 | 0.0161 | | | | | | | 6 8 | 0.0167 |
| 3 4 | 0.1376 | | | | | | | 6 9 | 0.0087 |
| 3 5 | 0.0352 | | | | | | | 6 10 | 0.0066 |
| 3 6 | 0.0225 | | | | | | | 7 8 | 0.0750 |
| 4 5 | 0.1399 | | | | | | | 7 9 | 0.0177 |
| 4 6 | 0.0416 | | | | | | | 7 10 | 0.0106 |
| 5 6 | 0.1589 | | | | | | | 8 9 | 0.0769 |
| | | | | | | | | 8 10 | 0.0216 |
| | | | | | | | | 9 10 | 0.0880 |

Table 2. *A-optimal design measures for block size* $k = 3$

| $v = 4$ | | $v = 7$ | | $v = 8$ | | $v = 9$ | | $v = 10$ | |
|---|---|---|---|---|---|---|---|---|---|
| $\tilde{\phi} = 8.5981$ | | $\tilde{\phi} = 31.6759$ | | $\tilde{\phi} = 42.6698$ | | $\tilde{\phi} = 55.2872$ | | $\tilde{\phi} = 69.5590$ | |
| Block | Mass | Block | Mass | Block | Mass | Block | Mass | Block | Mass |
| 1 2 3 | 0.3840 | 1 2 3 | 0.1860 | 1 2 3 | 0.1551 | 1 2 3 | 0.1329 | 1 2 3 | 0.1167 |
| 1 2 4 | 0.1160 | 1 2 4 | 0.0474 | 1 2 4 | 0.0289 | 1 2 4 | 0.0244 | 1 2 4 | 0.0186 |
| 1 3 4 | 0.1160 | 1 2 5 | 0.0404 | 1 2 5 | 0.0329 | 1 2 5 | 0.0261 | 1 2 5 | 0.0227 |
| 2 3 4 | 0.3840 | 1 2 6 | 0.0238 | 1 2 6 | 0.0343 | 1 2 6 | 0.0259 | 1 2 6 | 0.0207 |
| | | 1 2 7 | 0.0138 | 1 2 7 | 0.0126 | 1 2 7 | 0.0183 | 1 2 7 | 0.0149 |
| $v = 5$ | | 1 6 7 | 0.0138 | 1 2 8 | 0.0080 | 1 2 8 | 0.0093 | 1 2 8 | 0.0147 |
| $\tilde{\phi} = 14.6063$ | | 2 3 4 | 0.1119 | 1 7 8 | 0.0080 | 1 2 9 | 0.0048 | 1 2 9 | 0.0060 |
| Block | Mass | 2 6 7 | 0.0238 | 2 3 4 | 0.1044 | 1 8 9 | 0.0048 | 1 2 10 | 0.0028 |
| 1 2 3 | 0.3190 | 3 4 5 | 0.1534 | 2 7 8 | 0.0126 | 2 3 4 | 0.0901 | 1 9 10 | 0.0028 |
| 1 2 4 | 0.0831 | 3 6 7 | 0.0404 | 3 4 5 | 0.1238 | 2 8 9 | 0.0093 | 2 3 4 | 0.0821 |
| 1 2 5 | 0.0411 | 4 5 6 | 0.1119 | 3 7 8 | 0.0343 | 3 4 5 | 0.1054 | 2 9 10 | 0.0060 |
| 1 4 5 | 0.0411 | 4 6 7 | 0.0474 | 4 5 6 | 0.1238 | 3 4 7 | 0.0106 | 3 4 5 | 0.0922 |
| 2 3 4 | 0.1136 | 5 6 7 | 0.1860 | 4 7 8 | 0.0329 | 3 6 7 | 0.0106 | 3 4 7 | 0.0067 |
| 2 4 5 | 0.0831 | | | 5 6 7 | 0.1044 | 3 8 9 | 0.0183 | 3 4 8 | 0.0036 |
| 3 4 5 | 0.3190 | | | 5 7 8 | 0.0289 | 4 5 6 | 0.1044 | 3 6 7 | 0.0069 |
| | | | | 6 7 8 | 0.1551 | 4 8 9 | 0.0259 | 3 7 8 | 0.0036 |
| $v = 6$ | | | | | | 5 6 7 | 0.1054 | 3 9 10 | 0.0147 |
| $\tilde{\phi} = 22.3476$ | | | | | | 5 8 9 | 0.0261 | 4 5 6 | 0.0914 |
| Block | Mass | | | | | 6 7 8 | 0.0901 | 4 5 8 | 0.0069 |
| 1 2 3 | 0.2338 | | | | | 6 8 9 | 0.0244 | 4 7 8 | 0.0067 |
| 1 2 4 | 0.0612 | | | | | 7 8 9 | 0.1329 | 4 9 10 | 0.0149 |
| 1 2 5 | 0.0351 | | | | | | | 5 6 7 | 0.0914 |
| 1 2 6 | 0.0277 | | | | | | | 5 9 10 | 0.0207 |
| 1 5 6 | 0.0277 | | | | | | | 6 7 8 | 0.0922 |
| 2 3 4 | 0.1422 | | | | | | | 6 9 10 | 0.0227 |
| 2 5 6 | 0.0351 | | | | | | | 7 8 9 | 0.0821 |
| 3 4 5 | 0.1422 | | | | | | | 7 9 10 | 0.0186 |
| 3 5 6 | 0.0612 | | | | | | | 8 9 10 | 0.1167 |
| 4 5 6 | 0.2338 | | | | | | | | |

Table 3. *A-optimal design measures for block size* $k = 4$

| $v = 5$ | | $v = 7$ | | $v = 8$ | | $v = 9$ | | $v = 10$ | |
|---|---|---|---|---|---|---|---|---|---|
| $\tilde{\phi} = 10.2901$ | | $\tilde{\phi} = 22.0014$ | | $\tilde{\phi} = 29.5602$ | | $\tilde{\phi} = 38.2044$ | | $\tilde{\phi} = 47.9778$ | |
| Block | Mass | Block | Mass | Block | Mass | Block | Mass | Block | Mass |
| 1 2 3 4 | 0.3641 | 1 2 3 4 | 0.2614 | 1 2 3 4 | 0.2107 | 1 2 3 4 | 0.1846 | 1 2 3 4 | 0.1564 |
| 1 2 3 5 | 0.0757 | 1 2 3 5 | 0.0661 | 1 2 3 5 | 0.0332 | 1 2 3 5 | 0.0188 | 1 2 3 5 | 0.0093 |
| 1 2 4 5 | 0.1204 | 1 2 5 6 | 0.0570 | 1 2 3 6 | 0.0352 | 1 2 3 6 | 0.0012 | 1 2 3 6 | 0.0048 |
| 1 3 4 5 | 0.0757 | 1 2 6 7 | 0.0648 | 1 2 5 6 | 0.0125 | 1 2 3 7 | 0.0320 | 1 2 3 7 | 0.0237 |
| 2 3 4 5 | 0.3641 | 2 3 4 5 | 0.0831 | 1 2 6 7 | 0.0127 | 1 2 5 6 | 0.0376 | 1 2 3 8 | 0.0232 |
| | | 2 3 6 7 | 0.0570 | 1 2 7 8 | 0.0805 | 1 2 7 8 | 0.0176 | 1 2 5 6 | 0.0337 |
| $v = 6$ | | 3 4 5 6 | 0.0831 | 2 3 4 5 | 0.1168 | 1 2 8 9 | 0.0494 | 1 2 8 9 | 0.0105 |
| $\tilde{\phi} = 15.5521$ | | 3 5 6 7 | 0.0661 | 2 3 7 8 | 0.0127 | 2 3 4 5 | 0.0908 | 1 2 9 10 | 0.0416 |
| Block | Mass | 4 5 6 7 | 0.2614 | 3 4 5 6 | 0.0773 | 2 3 8 9 | 0.0176 | 2 3 4 5 | 0.0888 |
| 1 2 3 4 | 0.2975 | | | 3 4 7 8 | 0.0125 | 3 4 5 6 | 0.0927 | 2 3 9 10 | 0.0105 |
| 1 2 5 6 | 0.2364 | | | 3 6 7 8 | 0.0352 | 3 7 8 9 | 0.0320 | 3 4 5 6 | 0.0799 |
| 2 3 4 5 | 0.1686 | | | 4 5 6 7 | 0.1168 | 4 5 6 7 | 0.0927 | 3 8 9 10 | 0.0232 |
| 3 4 5 6 | 0.2975 | | | 4 6 7 8 | 0.0332 | 4 5 8 9 | 0.0376 | 4 5 6 7 | 0.0978 |
| | | | | 5 6 7 8 | 0.2107 | 4 7 8 9 | 0.0012 | 4 8 9 10 | 0.0237 |
| | | | | | | 5 6 7 8 | 0.0908 | 5 6 7 8 | 0.0799 |
| | | | | | | 5 7 8 9 | 0.0188 | 5 6 9 10 | 0.0337 |
| | | | | | | 6 7 8 9 | 0.1846 | 5 8 9 10 | 0.0048 |
| | | | | | | | | 6 7 8 9 | 0.0888 |
| | | | | | | | | 6 8 9 10 | 0.0093 |
| | | | | | | | | 7 8 9 10 | 0.1564 |

Table 4. *A-optimal design measures for block size* $k = 5$

| $v = 6$ | | $v = 8$ | | $v = 9$ | | $v = 10$ | |
|---|---|---|---|---|---|---|---|
| $\tilde{\phi} = 12.1358$ | | $\tilde{\phi} = 22.9128$ | | $\tilde{\phi} = 29.6043$ | | $\tilde{\phi} = 37.1256$ | |
| Block | Mass | Block | Mass | Block | Mass | Block | Mass |
| 1 2 3 4 5 | 0.3508 | 1 2 3 4 5 | 0.2441 | 1 2 3 4 5 | 0.2183 | 1 2 3 4 5 | 0.1958 |
| 1 2 3 4 6 | 0.0572 | 1 2 3 4 6 | 0.0031 | 1 2 3 4 6 | 0.0352 | 1 2 3 4 6 | 0.0181 |
| 1 2 3 5 6 | 0.0920 | 1 2 3 6 7 | 0.0593 | 1 2 3 6 7 | 0.0555 | 1 2 3 4 7 | 0.0164 |
| 1 2 4 5 6 | 0.0920 | 1 2 3 7 8 | 0.0910 | 1 2 3 7 8 | 0.0084 | 1 2 3 6 7 | 0.0055 |
| 1 3 4 5 6 | 0.0572 | 1 2 6 7 8 | 0.0910 | 1 2 3 8 9 | 0.0408 | 1 2 3 7 8 | 0.0273 |
| 2 3 4 5 6 | 0.3508 | 2 3 4 5 6 | 0.1025 | 1 2 6 7 8 | 0.0283 | 1 2 3 9 10 | 0.0547 |
| | | 2 3 6 7 8 | 0.0593 | 1 2 7 8 9 | 0.0408 | 1 2 7 8 9 | 0.0134 |
| $v = 7$ | | 3 4 5 6 7 | 0.1025 | 2 3 4 5 6 | 0.0878 | 1 2 8 9 10 | 0.0547 |
| $\tilde{\phi} = 17.1113$ | | 3 5 6 7 8 | 0.0031 | 2 3 4 7 8 | 0.0012 | 2 3 4 5 6 | 0.1061 |
| Block | Mass | 4 5 6 7 8 | 0.2441 | 2 3 4 8 9 | 0.0283 | 2 3 4 9 10 | 0.0134 |
| 1 2 3 4 5 | 0.2429 | | | 2 3 6 7 8 | 0.0012 | 3 4 5 6 7 | 0.0627 |
| 1 2 3 4 6 | 0.0115 | | | 2 3 7 8 9 | 0.0084 | 3 4 8 9 10 | 0.0273 |
| 1 2 3 6 7 | 0.1226 | | | 3 4 5 6 7 | 0.0370 | 4 5 6 7 8 | 0.0627 |
| 1 2 5 6 7 | 0.1226 | | | 3 4 7 8 9 | 0.0555 | 4 5 8 9 10 | 0.0055 |
| 2 3 4 5 6 | 0.1336 | | | 4 5 6 7 8 | 0.0878 | 4 7 8 9 10 | 0.0164 |
| 2 4 5 6 7 | 0.0115 | | | 4 6 7 8 9 | 0.0352 | 5 6 7 8 9 | 0.1061 |
| 3 4 5 6 7 | 0.2429 | | | 5 6 7 8 9 | 0.2183 | 5 7 8 9 10 | 0.0181 |
| 1 2 2 3 4 | 0.0562 | | | 1 2 2 3 4 | 0.0060 | 6 7 8 9 10 | 0.1958 |
| 4 5 6 6 7 | 0.0562 | | | 6 7 8 8 9 | 0.0060 | | |

11